\documentclass[a4paper,10pt]{article}
\usepackage[T2A]{fontenc}
\usepackage[utf8]{inputenc}
\usepackage[english]{babel}
\usepackage{amssymb}
\usepackage{amsmath}
\usepackage{umnbib}

\def\be{\begin{equation}}
\def\ee{\end{equation}}

\title{Variational formula in the Schottky model\\ of Riemann surfaces}
\author{\copyright A.B.Bogatyrev\thanks {Supported by RSCF, project  25-11-00144 
and INM RAS division of  MCFAM (Agreement with Russian Ministry of Science and Education 075-15-2025-347)} 
(INM RAS, Moscow State U., MCFAM, HSE U.)}
\date{}

\begin{document}
\maketitle

\begin{abstract}
A simple and efficient formula is proposed for varying Abelian integrals (including their periods) under variation of the generators of the classical Schottky group representing a Riemann surface.
\end{abstract}

Riemann surfaces are often used to solve various problems in the natural sciences, including mathematics, physics, and engineering. The question inevitably arises of efficiently computing answers given in terms of function-theoretic objects on surfaces: meromorphic functions, Abelian integrals, differentials (including quadratic ones), spinors, and so on. The Schottky model provides such an opportunity and is an alternative to using theta or sigma functions of several variables.

Suppose there are $2g$, $g=1,2,\dots$ disjoint disks $D_k$, $D'_k$, $k=1,\dots,g$ in the complex plane. We fix Möbius transformations $S_k\in~PSL_2(\mathbb{C})$, $k=1,\dots,g$, mapping the interior of $D'_k$ to the exterior of $D_k$. Such transformations are not uniquely determined by the discs, but they impose significant constraints both on the motions  $S_k$ themselves and on their interaction. For example, any motion $S_k$ is loxodromic with an attractive fixed point lying in the circle $D_k$ and the repelling point lying in $D'_k$). The motions $S_k$, $k=1,\dots,g$, freely generate the group $\mathfrak{S}$, called the classical Schottky group \cite{Sch87}. This group acts discontinuously on the Riemann sphere outside its limit set $\Lambda(\mathfrak{S})$, known as ``Cantor dust'' and having a fractal nature. The exterior of all $2g$ circles (with attached part of boundary) is the fundamental domain $\sf{R}(\mathfrak{S})$ of this action, and the manifold of orbits is the Riemann surface of genus $g$ obtained by gluing together the boundary circles of the fundamental domain using the generators of the group. We call thus obtained surface the Schottky model.

The function theory on a Riemann surface in the Schottky model can be efficiently constructed using linear Poincaré series (= automorphic forms of weight -2) \cite{Poi82}, representing Abelian differentials on the orbit manifold of the group $\mathfrak{S}$. For example, the Abelian differential of the third kind $d\eta_{zz'}$, invariant under the group action, with simple poles at the points $z,z'$ of the discontinuity domain $\mathbb{CP}^1\setminus \Lambda(\mathfrak{S})$, is obtained from the simplest differential with the same singularities on the sphere by summing over the action of all motions $S$ from the Schottky group:
\be
\label{detazz}
d\eta_{zz'}(u):=
\sum\limits_{S\in\mathfrak{S}}
\left\{
\frac1{Su-z}-\frac1{Su-z'}
\right\}
dS(u)
=
\sum\limits_{S\in\mathfrak{S}}
\left\{
\frac1{u-Sz}-\frac1{u-Sz'}
\right\}
du.
\ee
Such a series converges absolutely if the Schottky criterion \cite{Sch87} on the fundamental domain of the group is satisfied, for example, if the centers of all disks $D_k$, $D'_k$ are real. Holomorphic differentials are obtained from this construction by placing the poles in the same orbit of the Schottky group, and differentials with multiple poles are obtained by differentiating with respect to the pole position. Term-by-term integration of these series followed by exponentiation yields Schottky functions
\cite{Sch87,Bbook} from which it is easy to construct meromorphic functions on the orbit manifold. An alternative to this approach is the construction of the Klein's prime form $E(x,y)$ \cite{Mum83}, which also has a constructive representation in the Schottky model \cite{B09} and serves as the basis for the effective construction of function theory on a Riemann surface.

For applications, it is often necessary to find surfaces of a special type on which certain relations between distinguished functions, differentials, and so on hold. To do this, we must solve equations for the moduli of the surface that are parameterized e.g. by the generators of the Schottky group $\mathfrak{S}$.
Numerically, this can be done using Newton's method, which requires knowledge of the directional derivatives of the Schottky functions being calculated. Such formulas for derivatives, which allow for efficient implementation, are given below.

Let $\mathfrak{S}$ be the classical Schottky group with generators $S_1$, $S_2$, \dots, $S_g$,
represented in matrix form:
\be
\label{matrix}
S_j(u):=\frac{c_{11}u+c_{12}}{c_{21}u+c_{22}} \longrightarrow
\hat S_j:=\left(
\begin{array}{cc}
c_{11}~c_{12}\\c_{21}~c_{22}\\
\end{array}\right)
\in GL_2(\mathbb{C})
\ee

{\bf Theorem}
~~The following variational formulas are valid for Abelian integrals within given limits
\be
\label{main}
\delta\int\limits_z^{z'}d\eta= 
(2\pi i)^{-1}
\sum\limits_{l=1}^g \int\limits_{\partial D_l}
d\eta(u)d\eta_{zz'}(u) {\rm tr}[{\cal M}(u)\cdot\delta\hat S_l\cdot\hat S_l^{-1}]/du+o, 
\ee
Here $d\eta(u)$ is the Abelian differential (with any singularities) on the surface,
normalized over the cycles $\partial D_k$; the points $z,z'$ are from the fundamental domain, and
${\cal M}(u):=
\begin{array}{||cc||}
-u&u^2\\
-1&u\\
\end{array}
$; all intergrations are in the clockwise direction and the residual term $o:=O(\sum\limits_{sj}(\delta c_{sj})^2)$.

{\bf Remarks}
1) Variational formulas for the \emph{periods of Abelian differentials}
are obtained by placing the endpoints of the integration path in one orbit of the group and replacing
the differential $d\eta_{zz'}$ on the right-hand side of the formula with the corresponding holomorphic differential.\\
2) Similar variational formulas also exist for the Kleinian prime form in the Schottky model \cite{B09}.\\
3) For the Schottky group representing real hyperelliptic curves,
these formulas apparently first appeared in \cite{B99}; see also \cite{B03, Bbook}.
Here we establish their validity for an arbitrary Schottky group, including those with divergent linear Poincaré series of the type \eqref{detazz}.\\
4) The structure of the variational formulas \eqref{main} resembles that for the Hadamard variational formulas \cite{Had} for Green's functions. This  is not a surprise in view of the connection of the latter with Abelian integrals. The matrix ${\cal M}(u)$ was introduced by D. Hejhal \cite{Hej76} to represent the monodromy variation of projective structures, to which our formula \eqref{main} is also related.\\
5) The numerical efficiency of the variational formulas \eqref{main} is based on the work of D. Hejhal
\cite{Hej78}, in which integrals over the boundary components of the fundamental domain of the product of the inverse differential ${\cal M}(u)/du$ and the quadratic differential represented as a quadratic Poincaré series were calculated in finite form. Numerous examples of such calculations for real hyperelliptic curves are given in the book \cite{Bbook}. They are also extended to the more general case.

{\bf Sketch of the proof}. We consider the most simple case of variations for elements of period matrix.  It illustrates the idea of the proof while technical details for the general case may be easily recovered. 
Let $d\zeta_s$, $s=1,\dots,g$, be the basis of holomorphic differentials on the surface normalized as 
$\int_{a_s}d\zeta_j=\delta_{js}$, $a_s:=\partial D_s$. 
The variation of the generators $S_l$ we represent in the matrix form \eqref{matrix}:
$\hat{S_l}^*:=\hat{S_l}+\delta\hat{S_l}$. Consider the smooth mapping of the fundamental domain of the unperturbed group $\mathfrak{S}$ to the domain for the perturbed group $\mathfrak{S}^*$ which respects the boundary identifications:
\be
\label{qcmap}
f(u)=u + \sum_{l=1}^g\rho^l(u)[S_l^*\circ S_l^{-1}(u)-u]
\ee
where $\rho^l(u)$ is a smooth function equal to $1$ in a small vicinity of $D_l$
and vanishing in a slightly larger vicinity of the disc. This is a smooth map for the 
small enough perturbations of generators of the Schottky group $\mathfrak S$. Let us  compute its 
Beltrami coefficient keeping only linear terms in the perturbation expansion:
\be
\label{muPP}
\mu[f]:=f_{\bar{u}}/f_{u}=\sum_{l=1}^g(\rho^l_{\bar{u}}(u)(S_l^*\circ S_l^{-1}(u)-u)+o=
-\sum_{l=1}^g(\rho^l_{\bar{u}}(u)tr[{\cal M}(u)\cdot\delta\hat S_l\cdot\hat S_l^{-1}]+o.
\ee
Now we use Rauch variational formula \cite{Rauch} for the periods of holomorphic differentials along 
cycles $b_s$ which extend the set of cycles $a_s:=\partial D_s$, $s=1,\dots,g$, to 
a symplectic basis in 1-homologies on the surface. 
\be
\delta \int_{b_s}d\zeta_j= \frac1{2\pi i} ~\int_{\sf {R}(\mathfrak{S})} 
\frac{d\zeta_j(u)}{du}\frac{d\zeta_s(u)}{du}\mu[f](u) du\wedge \bar{du}+O(||\mu[f]||^2).
\ee
Inserting the principal part \eqref{muPP} of the Beltrami coefficient to the last formula and integrating by parts we get exactly formula \eqref{main} with $d\eta=d\zeta_j$, $d\eta_{zz'}=d\zeta_s$.

\end{document}